# Sub-Optimal Multi-Phase Path Planning: A Method for Solving Rubik's Revenge

Jared Weed Department of Mathematical Sciences
Worcester Polytechnic Institute
Worcester, Massachusetts 01609-2280
Email: jmweed@wpi.edu

*Abstract*—Rubik's Revenge, a $4 \times 4 \times 4$ variant of the Rubik's puzzles, remains to date as an *unsolved* puzzle. That is to say, we do not have a method or successful categorization to optimally solve every one of its approximately $7.401 \times 10^{45}$ possible configurations. Rubik's Cube, Rubik's Revenge's predecessor ($3 \times 3 \times 3$), with its approximately $4.33 \times 10^{19}$ possible configurations, has only recently been completely solved by Rokicki et. al, further finding that any configuration requires no more than 20 moves. [8] With the sheer dimension of Rubik's Revenge and its total configuration space, a brute-force method of finding all optimal solutions would be in vain. Similar to the methods used by Rokicki et. al on Rubik's Cube, in this paper we develop a method for solving arbitrary configurations of Rubik's Revenge in phases, using a combination of a powerful algorithm known as IDA* and a useful definition of distance in the cube space. While time-series results were not successfully gathered, it will be shown that this method far outweighs current human-solving methods and can be used to determine loose upper bounds for the cube space. Discussion will suggest that this method can also be applied to other puzzles with the proper transformations.

## I. INTRODUCTION

Rubik's puzzles have been an interest to children and adults alike for many decades since the original Rubik's Cube invention in 1974 by Erno Rubik. Today, many variations of Rubik's Cube exist and each provide their own challenges when solving. Interestingly, such a large following for these puzzles have developed over the years, that entire competitions revolve around the speed at which the puzzles can be solved, using hundreds of memorized algorithms designed for speed and recognition. The most popular of which have been the Rubik's Pocket, Rubik's Cube, Rubik's Revenge, and Rubik's Professor. The World Cube Association, since 2003, has conducted competitions and kept record of results.

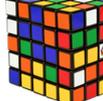

Fig. 1. Current World Records in 2015

While many enjoy solving these puzzles by hand using a set of well-established algorithms, mathematicians find interest in understanding the properties the puzzles have on a more fundamental level. Since each puzzle is fully observable, we can derive useful results based on the theory alone. As an example, the algorithms used in speed-solving by humans are terribly suboptimal, in that the necessary amount of moves to correctly solve the scrambled configuration is much less than the required moves by the algorithms. This is by virtue of each algorithm being generalized and not specific to one particular configuration. This however begs the question: What if we were gods? What if we had the ability to *know* how to solve any configuration of a puzzle in the least amount of moves possible, through the use of some clever algorithm? This hypothetical notion is defined as **God's Algorithm**, one that unfortunately we as humans do not know.

Throughout this paper, we will be adopting World Cube Association's (WCA) Single-Turn notation to annotate how twists are executed on Rubik's Revenge. A quarter turn in the clockwise direction of the outer left, right, front, back, up, or down faces is aptly noted as $L, R, F, B, U, D$ respectively. Similarly, a quarter turn in the clockwise direction of the inner faces are noted as $l, r, f, b, u, d$. A turn that is instead in the counter-clockwise direction is represented with an additional tick (such as the move $R'$ - a counter-clockwise turn of the outer right face). Lastly, a half turn (as opposed to quarter turn) is represented with an additional 2 (such as the move $u2$ - a half-turn of the inner up face).

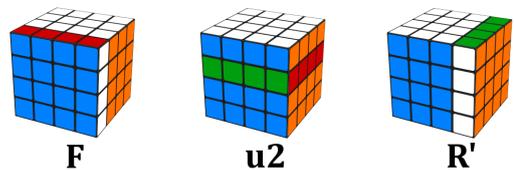

Fig. 2. Three Independent Moves Using WCA Notation

## II. RELATED WORK

Much of the method employed here stems from previous work developed on Rubik's Cube. The development that started it all occurred in 1981 when Morwen Thistlethwaite proposed in a paper that the Rubik's Cube can be solved using four independent phases. [9] The idea was that the cube, given any scramble, could be maneuvered into a configuration using quarter turns of all the faces, that could then be maneuvered into a configuration using half turns of some of the faces, and so on until the cube was solved. The choice of maneuvers was not arbitrary: He noted that as each phase was solved, the groups of available twists of the cube were nested in one another. For the first phase, the moves $L, R, F, B, U, D$ are permitted, and using these moves, the cube should be maneuvered into a configuration such that it can be solved using the moves $L, R, F, B, U2, D2$. The following figure gives a documentation of Thistletwaite's process. Note here that each subgroup denoted by $G_i$ are the moves available to the given phase, as mentioned earlier. It is obvious that these subgroups are indeed nested such that $G_i \supset G_{i+1}$. The number of configurations immediately following the subgroups are the total configurations that are possibly generated by that subgroup. Clearly, $G_1$ generates the entire cube space, which is known to have approximately $4.33 \times 10^{19}$ valid configurations. Following this is the number of coset elements. While this terminology will be explored later in the paper, we may naively define this as the number of unique classes of elements existing in $G_i$ not currently solvable in $G_{i+1}$. Lastly, a diagram is provided showing the state of a particular scramble as it is maneuvered through the phases of this algorithm.

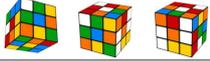

Fig. 3. Thistlethwaite's Four-Phase Process

In 1985, Richard Korf proposed a new method of tree searching called *Depth-first Iterative-deepening*, which could be used with a heuristic method such as A* (which he aptly coined IDA*). [3] Using this method, in 1992 Herbert Kociemba expanded on the theory of Thistlethwaite's four-phase process for Rubik's Cube and simplified it to only two-phases. [2] This was possible by combining the first two phases of Thistlethwaite's process into one phase, and also the last two-phases into one phase. From the previous figure, it would then be clear that the number of "coset elements" drastically increases in this new method. Note, however, that using IDA* on this bi-part search space allowed for considerably real-time results, and found that solutions typically of 45 to 52 moves in Thistlethwaite's process were reduced to solutions on average of 26 moves in this method.

Several years later, Korf showed that optimal solutions could be found for Rubik's Cube with the help of pattern databases, while also employing his previously developed IDA* algorithm. [4] The pattern databases, in this case, are related to the coset elements mentioned previously, in that each of these elements can have their distances computed and stored into a large table, each with a unique identifier. [1] Then, for any scrambled cube, it is identified by certain properties and its distance quickly found by a lookup in the table. This use of lookup tables drastically improved computation time when searching the cube space in conjunction with using iterative-deepening. Korf was able to show that a pattern database along with an efficient search algorithm provided him with optimum solutions as long as 18 moves.

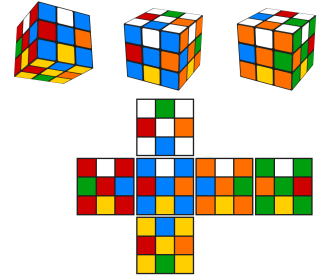

Fig. 4. The "Superflip": First Known Twenty-Move Configuration

Lastly, over the course of the years 2009 to 2013, a culmination of efforts from Korf, Kociemba, Rokicki, and others, found the optimal solution to every valid configuration of Rubik's Cube and determined that the maximum number of moves necessary to solve one is 20. [7][8] While it is clear that the computation necessary to do so required a timespan of several years, the method by which these solutions were found came largely from the development of methods by Korf using IDA* and the pattern databases (which Kociemba incidentally adopted for his two-phase solver in recent years).

## III. PROBLEM STATEMENT

It is clear that Rubik's Cube has, by definition, been solved and categorized for every valid configuration. However, we note that this success came only just a couple of years ago. Rubik's Revenge, its successor in size, has yet to be solved, and this is primarily due to the sheer difference in size of the configuration space when comparing it to that of Rubik's Cube. With minimal effort, it can be shown that the approach used by Rokicki et. al to solve Rubik's Cube cannot be applied directly to the $7.401 \times 10^{45}$ configurations, as a single search tree would require a depth upwards of 34 and a breadth of nearly 24. Such an expanse would require a near-infinite amount of time in the current computing power available.

Instead, in this paper we investigate the following question: Can a series of sub-problems be planned optimally for the Rubik's Revenge to sub-optimally solve any valid configuration? This means that we will not be searching for God's Algorithm exactly, but rather cheat the definition and find generalized versions of the algorithm for each sub-problem. We note that this approach is similar to that of previous works of Thistlethwaite and Kociemba, but will necessarily require the approach of Korf and Rokicki when applying IDA* to find solutions at greater depths. In the following section, we will define several components to this method. Firstly, we will develop some basic theory about the puzzle which will aid us greatly in sub-dividing the problem. Next, a definition of a multi-phase method will be given, followed by a definition of the metric used to have a notion of "distance from solved" for each phase. Finally, pseudo-code for the two-part algorithm of IDA* will be provided, showing how the method approaches deriving a solution for any configuration.

## IV. METHODOLOGY

*Theory*

As mentioned in the sections above, there are approximately $7.401 \times 10^{45}$ possible valid configurations of Rubik's Revenge. [5] There are eight corner pieces, which are unrestricted in their orbit, resulting in 8! possible permutations of the corners. Furthermore, each corner piece has one of three possible rotations, with one such piece having its rotation defined by the other seven (the parity of the corners must always be even), resulting in another factor of $3^7$. Thus in total, the corner pieces have $8! \cdot 3^7$ unique configurations. Additionally, there are twnety-four edge pieces, which are also unrestricted in their orbit, resulting in 24! possible permutations of the edges. The edge pieces, however, cannot be flipped singly since the internal structure of Rubik's Revenge is asymmetrical for each edge pair.

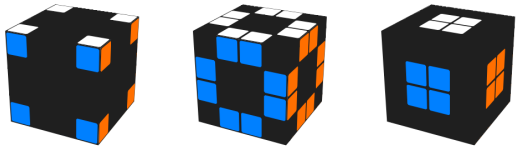

Fig. 5. The corners, edges, and centers of Rubik's Revenge, respectively

For the twenty-four center pieces, which also have 24! possible permutations. It is important to note that the center pieces on a solved Rubik's Revenge are indistinguishable; that is, there is no designation where each of the four center pieces per face belong. Thus we can reduce the total possible permutations by a factor of 4! for each face, resulting in a total reduction factor of $(4!)^6$. Lastly, we note that the orientation of the cube is not fixed in space. Unlike Rubik's Cube, which has fixed centers to designate the color scheme, Rubik's Revenge has no such fixed center, thus we can additionally reduce the permutations by a factor of 24 to account for symmetry.

When putting the above calculations together, we find that the total number of possible configurations is

$$\frac{8! \cdot 3^7 \cdot 24! \cdot 24!}{(4!)^6 \cdot 24} = \frac{8! \cdot 3^7 \cdot 24!}{24} \approx 7.401 \times 10^{45}$$

as shown previously.

Given that there are 12 independent slices for Rubik's Revenge, each of which can be twisted a total of 3 different ways, we have a total of 36 unique twists that can be executed on the cube. However, if we consider a *sequence* of twists exacted on the cube, we find that there are only 24 unique twists when performed in succession. This is not immediately obvious, however it can be easily understood. Consider the sequence of two moves $F\ F2$, which is a quarter-turn of the outer front face, followed by a half-turn of the outer front face. It comes as no surprise that this two-twist sequence can be identified by the single twist $F'$. Similarly, we can consider a different type of two-twist sequence, such as $R2\ L2$, which is a half-turn of the outer right face, followed by a half-turn of the outer left face. The resulting configuration is exactly identical to the resulting configuration if we swapped the order of the sequence, namely $L2\ R2$. In this case, the sequence is commutative, because the two slices being performed on are independent of one another. Therefore, for any arbitrarily chosen twist of an arbitrarily chosen slice, we can ensure the sequence is non-commutative or non-simplifiable (at least initially) by disallowing subsequent twists that occur on a parallel slice. Thus this results in the figure of 24 from before.

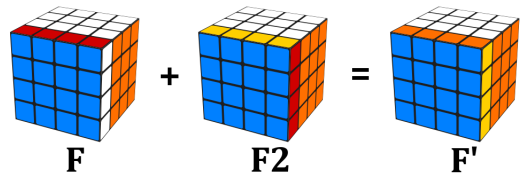

Fig. 6. The result of two twists being equivalent to one.

Using these facts, we can show by a fairly straightforward pigeonhole argument that there exists a configuration of the cube that requires at least 34 moves to solve. Since 36 twists are available for the first move in a sequence, and 24 for all subsequent twists, we are curious to solve the inequality

$$36 \cdot 24^{n-1} \geq 7.401 \times 10^{45}$$

where $n$ is the number of twists. A quick calculation shows that $n$ is necessarily larger than 34.

*Solving in Phases*

By the information above, it is obvious that a brute-force method involving one search tree is completely intractable. We note that the tree will be at least 34 nodes deep, and have an approximate branching factor of 24. This indicates that the method will require the problem be solved in phases. We do this by dividing the cube group into nested subgroups $G_0 \supset G_1 \supset \ldots \supset G_n$, where $G_0$ represents the cube group, and $G_n$ represents the identity element (the solved configuration). Each subgroup is then defined by the twist generators, which can generate a subset of the total possible configurations. Since the subgroups are nested, this implies that as we move from $G_i$ to $G_{i+1}$ we restrict the twist generators available to that class. As an example, recall Thistlethwaite's four-phase scheme from above, and note that the subgroups are

$$G_0 = \{L, R, F, B, U, D\}$$
$$G_1 = \{L, R, F, B, U2, D2\}$$
$$G_2 = \{L, R, F2, B2, U2, D2\}$$
$$G_3 = \{L2, R2, F2, B2, U2, D2\}$$
$$G_4 = \{I\}$$

from this, we can see that his four-phase scheme indeed has nested subgroups and that $G_0$ represents the twist generators of all possible Rubik's Cube configurations.

Next, we consider an arbitrary configuration that can be generated by the group $G_i$. Looking ahead to the following group $G_{i+1}$, we determine the unique characteristics about this configuration that cannot be solved using the generators of $G_{i+1}$. Interestingly enough, the elements contained in the right coset space $G_{i+1} \backslash G_i$ are all the possible unique classes of configurations that define these characteristics between the groups $G_i$ and $G_{i+1}$. This is not intuitive, and difficult to visualize given the groups we are operating in. The following is an example of right coset spaces:

*Let G be the group of positive integers modulo 8, i.e,*

$$G = \{0, 1, 2, 3, 4, 5, 6, 7\} = \mathbb{Z}_8^+$$

*and let $H = \{0, 4\}$, which is indeed a subgroup of G. Given that the operator on the integers is addition, the following elements are contained in the right coset group $H \backslash G$ :*

$$H + 0 = \{0, 4\}$$
$$H + 1 = \{1, 5\}$$
$$H + 2 = \{2, 6\}$$
$$H + 3 = \{3, 7\}$$

Using the example above, we can think of the elements in the right coset space $G_i + 1 \backslash G_i$ as classes of elements that identify particular characteristics of the cube configuration. For instance, the first right coset group of Thistlethwaite's scheme has 1024 elements, which correspond to the 1024 possible configurations the twelve edges can be in. Note that multiple configurations can have the same configuration of the twelve edges, and thus get grouped into the same class element. This is the essence of the right coset space.

There is then an identity element in each of these right coset spaces, which is a class of configurations that are indeed solvable in the group $G_{i+1}$. We can thus reduce the problem to planning a sequence of moves from the class element corresponding to the configuration in $G_i$ to the identity element in $G_{i+1} \backslash G_i$, which brings the cube to a configuration solvable in $G_{i+1}$. The planning necessary within each right coset space will them be developed for each phase.

*Metric Definition*

Prior to constructing an algorithm, we must first derive a definition of *distance* in the cube space, to understand how "far away" a certain configuration is from being solved in a phase. This is necessary since a planning method such as A* (as part of the IDA* algorithm) implements an additional heuristic to better plan towards a solution. We thus develop a notion of **twist distance**. For any piece of the cube, find the minimum number of twists necessary to correctly position and orient it into a given phase configuration. That is, if the configuration is currently operating in $G_i$, find the number of twists for that piece such that it is in a configuration from the identity element of $G_{i+1} \backslash G_i$. We then define **twist distance** as the sum over all pieces characteristic of that phase.

As a heuristic, this current definition of twist distance is not admissible. To be admissible, we must divide the sum by 16. If we consider an arbitrary twist of a slice, we either disturb 12 pieces (for an inner slice) or 16 pieces (for an outer slice). Since a distance heuristic is admissible if it underestimates the true distance for any configuration, we must divide by the greatest disturbance. We can improve this metric by instead calculating the twist distance for three separate sets: the corners, the edges, and the centers. since the maximum disturbance of each of these sets is 4, 8, and 8 respectively, each corresponding twist distance must be divided by this factor. We can then take the maximum of these three values as our measured *distance*, which will be used in the following algorithm.

*Algorithm*

We then introduce iterative-deepening A* (IDA*) as our algorithm of choice, following similar construction used by Kociemba and Korf on Rubik's Cube. We want to find an optimal solution for each phase, however depth-first search doesn't guarantee an optimal solution (or one at all) for an arbitrary depth, and breadth-first search, with a branching factor of 24, would require a vast amount of memory. Iterative-deepening combines the luxury of both types of search while also guaranteeing an optimal solution without the need for memory. This is done by simply iterating the depth for each phase, in a pseudo-breadth-first search fashion, attempting to find a solution at depth $d$. If one cannot be found, then simply search again for a solution at instead depth $d + 1$. We then combine this search with the A* planner (and the above heuristic) to determine the branching of successive twists: If a solution must exist at depth $d$ and we have already conducted $n$ twists, then we can disallow a twist if the resulting configuration has a twist distance that is greater than $d - n$. This action performs a much-needed pruning of the search space if it's clear a solution will not be found within the restricted depth.

Combining all of the information above results in the following two algorithms. The first algorithm, $RevengeSolver$, corresponds to the more general-level planning, wherein we iterate through each phase to find the optimal solution. For a given valid cube configuration $C$, we begin with an unfilled *solution* plan. Then for each phase, we slowly increase the $depth$ at which the phase operates (i.e., does a $depth = 3$ solution exist? If not, try to find a $depth = 4$ solution, etc.). If a $path$ is returned from $PhaseSearch$ (the second algorithm), then the phase has produced a solution. In this case, we transform the current configuration of the cube using the sequence of twists in $path$ and append the $path$ to the current solution sequence. At this point, the current configuration is now solvable in the next phase, so move into the next phase until no more exist. Once this has completed, we must have found our solution, so we return $solution$ which is the sequence of twists necessary to solve the configuration.

The second algorithm, $PhaseSearch$, does the actual phase planning, given a configuration $C$, the current $phase$, and the $depth$ of twists remaining to be explored. Note that this algorithm is recursive, so we first introduce a recursion check: if $depth = 0$ (meaning we have no more allowable twists to

**Algorithm 1:** RevengeSolver($C$)

**Data**: $C$, a valid cube configuration
**begin**
    $solution \longleftarrow$ Null
    **for** $phase = 0 \ldots N$ **do**
        **for** $depth = 0 \ldots \infty$ **do**
            $path \longleftarrow$ PhaseSearch($phase, C, depth$)
            **if** Success **then**
                Transform($C, path$)
                $solution \longleftarrow solution + path$
                **break**
    **return** $solution$

explore), we calculate the twist distance $C$ in the current phase from the identity element of that phase. If this distance is 0, then we have necessarily arrived at the solution for this phase. If, however, $depth$ is not 0, then we again calculate the twist distance of $C$ in the current phase from the identity element of that phase. If this distance is less than $depth$ (i.e., a solution is still feasible), we can apply all the available twists from that phase to the configuration, and recurse through $PhaseSearch$, this time reducing $depth$ by 1.

**Algorithm 2:** PhaseSearch($phase, C, depth$)

**Data**: $C$, a valid cube configuration
       $phase$, current operational phase
       $depth$, depth left to travel
**begin**
    **if** $depth = 0$ **then**
        **if** Distance($C, phase$) $= 0$ **then**
            **return** Success
    **else**
        **if** Distance($C, phase$) $\leq depth$ **then**
            **for** $twist \in$ Successor($phase$) **do**
                $C_{new} \longleftarrow$ Transform($C, twist$)
                PhaseSearch($phase, C_{new}, depth - 1$)

## V. Theoretical Results

Theoretically speaking, the number of phases implemented in the method is inversely correlated to plan time. This is caused by several latent features present in the algorithm. Firstly, by introducing more phases, we inadvertently reduce the cardinality of the right coset spaces (although the number of right coset spaces is exactly the number of phases). This implies that the search space within each right coset space is smaller, requiring fundamentally less computation time and search. However, since we are searching through more phases, this necessarily increases the average total moves necessary to solve the cube. While we can suspect that the average number of moves within each phase decreases, the number of phases places a larger role in the monotonicity of the upper bound.

Conversely, by reducing the number of phases, we greatly increase the size of the right coset spaces. Intuitively, this is because a reduction in phases necessarily requires that more cases are examined that must be solved before moving into the next phase. Thus a much greater time is expended searching the space for a solution. However, the benefit to reducing the number of phases is that, as mentioned about the monotonicity of the upper bound, less moves on average will be required to solve the cube. In fact, as the number of phases is reduced to 1, we converge on God's Algorithm. It should be noted however, that even a reduction by 1 phase exponentially increases the amount of computation time needed.

Lastly, of importance is the fact that by virtue of the IDA* algorithm, each phase is guaranteed to be solved optimally. This means that if a configuration that meets the necessary conditions of a particular phase was found to have an $N$-move solution, we can guarantee that no solution of $N - 1$ or less existed for that configuration. There are, however, caveats to this fact. Firstly, we cannot guarantee that only one optimum solution exists for any configuration within a particular phase. In fact, it is quite likely that several optimal solutions exist. This comes with an interesting effect on the method itself. Different optimal solutions may yield different final configurations when maneuvered into the next phase. This means that the solution to successive phases are affected by their predecessors. Having two optimal solutions for one particular phase does not necessarily imply that their respective solutions in the following phase are of the same length.

## VI. Example Problem

While Bruce Norskog had found after tedious computation that an upper bound can be achieved on the Rubik's Revenge with only 79 single-twist moves, other cubing enthusiasts such as Charles Tsai worked on developing phase-like solutions, borrowing from the theory that Thistlethwaite previously proposed. [6] One such solution structure is an eight-step algorithm, which is in fact a two-step four-phase solver. [10] The first four phases work as you'd expect: successive subgroups continue to restrict the allowable twists of all twelve slices:

$$G_0 = \{R, L, F, B, U, D, r, l, f, b, u, d\}$$
$$G_1 = \{R, L, F, B, U, D, r, l, f2, b2, u2, d2\}$$
$$G_2 = \{R2, L2, F, B, U, D, r2, l2, f2, b2, u2, d2\}$$
$$G_3 = \{R2, L2, F2, B2, U, D, r2, l2, f2, b2\}$$

Once the fourth phase is successful, the intriguing aspect of this Tsai's algorithm is that the Revenge is now in a configuration known as a **reduction**. A Revenge in reduction means that, if inner slice twists are completely restricted, then the Revenge is isomorphic to the 3x3 Rubik's Cube. If that is the case, then Thistlethwaite's four-phase algorithm can be applied, and this is indeed the second-step of the algorithm and the last four phases:

$$G_4 = \{R, L, F, B, U, D\}$$
$$G_5 = \{R2, L2, F, B, U, D\}$$
$$G_6 = \{R2, L2, F2, B2, U, D\}$$
$$G_7 = \{R2, L2, F2, B2, U2, D2\}$$

After completing the eight phases corresponding to the above subgroups, any Rubik's Revenge can be solved. A visual example is as follows. Let us first apply a 35-move WCA random scramble to a fully solved Rubik's Revenge.

$$F\ d2\ u\ B\ 'F'\ u\ U\ L2\ r2\ f'$$
$$R\ D\ U\ l2\ u'\ r2\ R'\ B'\ f'\ l2$$
$$b\ L\ b\ u'\ f\ l\ b\ d\ U2\ B\ D'\ d2$$

This yields the following initially scrambled configuration:

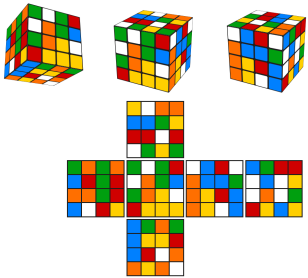

Fig. 7. An initially scrambled Revenge

*Phase 1*

In Phase 1, no twists are restricted when approaching $G_1$. To ensure the cube can be solved using the restrictions of $G_1$ the center pieces for both the right and left faces must contain only the two colors designated for those faces. This is because in $G_1$, there is no way to change the eight pieces set in the left and right centers. Phase 1 is solved as

$$B2\ u\ F2\ D\ r'\ f$$

From the initially scrambled configuration shown before, this results in the following configuration.

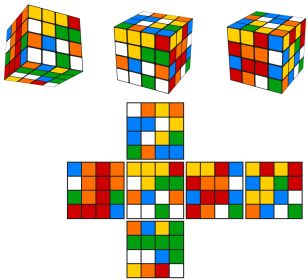

Fig. 8. After Phase 1: $G_1 = \{R, L, F, B, U, D, r, l, f2, b2, u2, d2\}$

*Phase 2*

Note that from the above figure, the left and right centers contain only red and orange pieces. In Phase 2, we now restrict nearly all the inner slices to half twists, with the exception of the $u$ and $d$ slices. To bring this configuration to one that can be solved using the restrictions of $G_2$, we must ensure that all edge pairs (currently unpaired) have the same parity. Furthermore, since $G_2$ restricts all inner slices to half twists, the center pieces for all opposing faces must contain only the two colors designated for those faces. Phase 2 is solved as

$$R2\ D2\ r\ B2\ L\ F'\ U'\ l$$

Which results in the next configuration

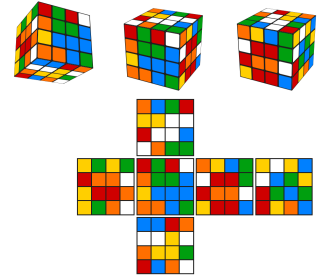

Fig. 9. After Phase 2: $G_2 = \{R2, L2, F, B, U, D, r2, l2, f2, b2, u2, d2\}$

*Phase 3*

In similar fashion as Phase 2, Phase 3 will now solve the centers for all opposing faces into an equivalent scheme made for the right and left centers. Furthermore, the configuration of each center face is such that they contain columns, which allows for the centers to be solved using only half-turn twists of the inner slices. To have the cube in a configuration that can be solved in $G_3$ we must also ensure that four arbitrary edges have been paired together and placed in the "middle" of the cube. Geometrically, these are the four edge positions not located on the top or bottom. The choice of these four pairings is unbounded so thus the shortest sequence of moves that does so is chosen. Phase 3 is solved as

$$R2\ B\ U'\ r2\ b2\ d2\ D'\ F\ r2\ B$$

Resulting in the following configuration, which you may note has four paired edges.

*Phase 4*

Phase 4 brings the configuration to reduction. Again, this means that it will then be in a configuration solvable using only the outer slices. This is done by pairing up the remaining eight edges and by solving each center so that the correct color is on each face. Furthermore, since $G_4$ will disallow any inner slice twists, we must also ensure that the parity of the corners and the edges match. Without this condition, the configuration will not match the parity of a Rubik's cube and thus not be isomorphic. This occurs since the centers are indistinguishable. The following moves solve Phase 4:

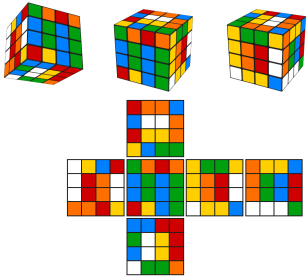

Fig. 10. After Phase 3: $G_3 = \{R2, L2, F2, B2, U, D, r2, l2, f2, b2\}$

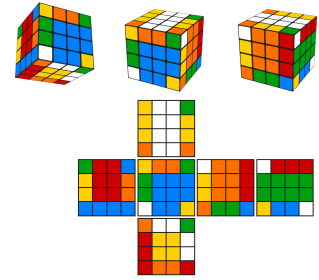

Fig. 12. After Phase 5: $G_5 = \{R2, L2, F, B, U, D\}$

$$F2\ l2\ U2\ R2\ B2\ U2\ r2\ U\ f2$$

It is clear at this point that all edges are paired and centers solved, resulting in a reduction.

$$U\ F\ D\ L2\ B'\ F'\ U\ B$$

Note that the "middle" consists of the colors Orange, Red, Blue, and Green, and thus the four edges corresponding to those colors are positioned in that layer.

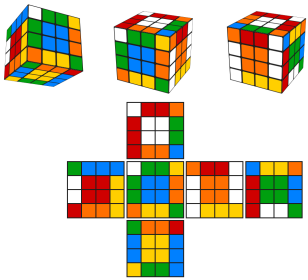

Fig. 11. After Phase 4: $G_4 = \{R, L, F, B, U, D\}$

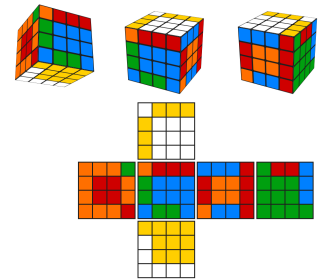

Fig. 13. After Phase 6: $G_6 = \{R2, L2, F2, B2, U, D\}$

*Phase 5*

Since we are now in reduction, we no longer consider twists of the inner slices. For this phase, we must first correct the parity of the edges, and we can do this without restriction of twists for the outer faces. If we consider $G_5$, we are restricted to half turns on the right and left faces. This means that the orientation of edges cannot be changed when considering its final position. Thus each edge must have a correct orientation before being solvable in $G_5$. Phase 5 is solved using the following moves:

$$L\ B\ L\ F'\ R$$

Unfortunately, this phase is not readily visualized, but note that each edge has "good" orientation.

*Phase 6*

In similar fashion to correcting parity of the edges in Phase 5, in Phase 6 we must ensure that we have the correct parity of the corners. Looking at $G_6$, we are further restricting the $F$ and $B$ faces to half turns as well, meaning that the corners cannot be reoriented when considering their final positions. Also, given this restriction, the four edges with a final position in the "middle" layer (this is the same as in Phase 3) must be placed in that layer. Phase 6 has the following solution:

*Phase 7*

In Phase 7, we place the remaining edges in their correct layers, without regard to the final position. Furthermore, since $G_7$ is a restriction of all faces to half-twists only, the permutation of the corners and edges must match: If the corners have even permutation, then so must the edges; if the corners have odd permutation, then the edges must also. We must satisfy this condition before a configuration can be solved in $G_7$ since half turns restrict us to even permutation twists only. With the following twists, Phase 7 is solved:

$$R2\ F2\ U\ F2\ U'\ L2\ U\ R2\ U$$

Visually, the cube is nearly solved, and we only need to solve opposing faces for the correct color.

*Phase 8*

Unsurprisingly, Phase 8 solves the cube, using only half-turns of the outer faces. Since we previously corrected for the permutation and orientation of both the edges and corners, we can fully solve the cube using the following sequence of moves:

$$R2\ D2\ R2\ U2\ R2\ B2\ D2\ R2\ F2\ L2\ B2\ D2$$

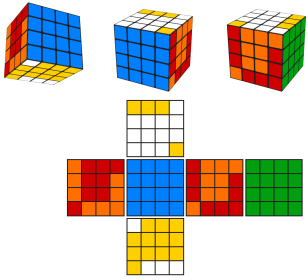

Fig. 14. After Phase 7: $G_7 = \{R2, L2, F2, B2, U2, D2\}$

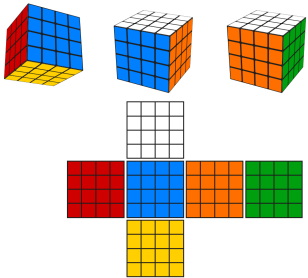

Fig. 15. Completion of Phase 8 results in solved Revenge

## VII. CONCLUSION

In this paper we investigated if a method can be developed to solve Rubik's Revenge in phases by considering the works done on Rubik's Cube. Indeed by confirmation of enthusiasts such as Bruce Norskog and Charles Tsai, phase solutions akin to Thistlethwaite's proposed scheme exist and are completely computational in real time. We have shown that an eight-phase solver exists, and that for the arbitrarily chosen scramble, only 67 moves are necessary to solve it using this scheme. While this is only one such example, many other phase solutions are bound to exist; an existence only limited by the computation time necessary for each phase given the properties outlined in the methodology section.

It is clear that a multi-phase solution is useful and perhaps most efficient when subdividing the problem. With the group-theoretic properties discussed in this paper, a brute-force solution for each possible configuration is impossible given our current computing power. It is then necessary to subdivide the problem into solvable sub-problems, which can then be solved in phases. The method employed has the benefit of providing real-time results that are sub-optimal by nature. Should a least-moves solution be necessary, this method can yield an upper bound on the number of moves constrained only by the construction of the subgroups.

Furthermore, given that the method discussed is derived within the properties of group theory, we may use it to an advantage to calculate worst-case analyses. As mentioned above, this method can yield an upper bound on the number of moves given a definition of the subgroups. Because of this, one may determine the strength of this construction by calculating this upper bound and comparing to other constructions. In this light it is possible to experiment with finding the most optimum N-phase solver.

Lastly, we note that this method can be applied to other puzzles of similar design. As we have seen, this method stems from its development for Rubik's Cube, and has been employed for Rubik's Revenge. Larger $N \times N$ cubes such as Professor's Cube share the same characteristics and group properties as already discovered and could quickly adopt the methods used for Rubik's Revenge. Furthermore, non-square cubes such as the Tetrahedron and Megaminx are suspect to also be able to adopt this method, as their properties are fully observable. Thus, similar constructions can be made.

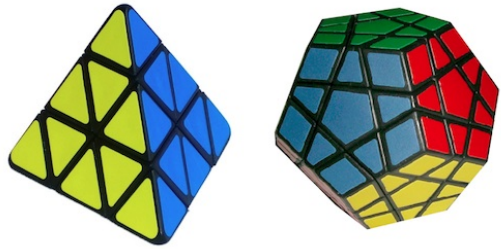

Fig. 16. Tetrahedron and Megaminx Puzzles

## VIII. FUTURE WORK

Given that the phases are independent of one another, one small improvement that can be made is to consider the boundaries between each phase to reduce redundancies. As an example, if Phase $N$ had the solution

$$U\ B\ R'\ L2\ U$$

and then Phase $N+1$ had the solution

$$U2\ R2\ L2\ F2$$

we may note that between the boundary, a $U$ twist followed by a $U2$ twist occurs, which is equivalent to $U'$. However, $U'$ may not necessarily be a move "allowed" in phase $N+1$. This is a correction that would be applied post-solution.

Secondly, the metric used in the method is such that all twists, disregarding the distance of rotation, have equal weighting. Other metrics can be investigated which may reveal different results. For instance, a penalty can be added to twists that are half-turns instead of quarter-turns. Alternatively, one can consider a block-turn metric, where both the outer and inner slices of a face can be turned together (as if they were glued), with the same weight as a single-slice turn. Such a metric has been used when constructing algorithms, where an example notation is $Ll$ or $(Uu)'$.

Since the orientation of Rubik's Revenge is not fixed in space (note that in this method an arbitrary orientation was given), another improvement that can be made is to allow rotations of the entire cube in the solution planning. Doing this can result in at least a reduction factor of 24, since we then can approach any of the 24 symmetries by virtue of not having a fixed orientation. Understandably, introducing rotation in the planning may result in a need to reconfigure the scheme for which the cube is identified, and adjusted so that it is rotation-dependent.

Lastly, note it was previously mentioned that the optimality of successive phases are dependent on the solution of their predecessors. It is possible to introduce an additional relaxation factor, which can potentially assist in finding more appropriate solutions over time. As an example, once a solution is found requiring a total of $t_1 + t_2 + t_3 + \ldots$ twists (where $t_i$ represents the number of twists to solve phase $i$) we can introduce a relaxation into the first phase, allowing a solution for that phase to require $t_1 + 1$ twists. By doing this, the resulting configuration may be solvable in the next phases requiring less moves than before.

ACKNOWLEDGMENT

For the many Rubik's puzzle diagrams, the author would like to thank David Whogg, for his Interactive Cube program written in Python. It is available for public use under the GNU General Public License. The program can be found on GitHub at the following link:

https://github.com/davidwhogg/MagicCube